\documentclass{gtart}

\def\ifplaintex{\expandafter\ifx\csname documentclass\endcsname\relax}

\ifplaintex 
\else
\fi

\def\gtm{{\mathsurround=0pt\it $\cal G\mskip-2mu$eometry \&\ 
$\cal T\!\!$opology $\cal M\mskip-1mu$onographs}}    

\def\gtp{{\mathsurround=0pt\it $\cal G\mskip-2mu$eometry \&\ 
$\cal T\!\!$opology $\cal P\!$ublications}}  

\def\recd{{\small Received:\qua\receiveddate\ifx\reviseddate\relax
\else\qquad Revised:\qua\reviseddate\fi\par}} 

\def\volumenumber#1{\def\thevolumenumber{#1}}
\def\volumeyear#1{\def\thevolumeyear{#1}}
\def\volumename#1{\def\thevolumename{#1}}
\def\papernumber#1{\def\thepapernumber{#1}}
\def\pagenumbers#1#2{\def\startpage{#1}\def\finishpage{#2}}
\def\published#1{\def\publishdate{#1}}
\def\received#1{\def\receiveddate{#1}}
\def\revised#1{\def\reviseddate{#1}}
\def\accepted#1{\def\accepteddate{#1}}

\def\asciiaddress#1{\def\theasciiaddress{#1}}

\let\\\par
\let\thevolumenumber\relax\let\thepapernumber\relax
\let\thevolumeyear\relax\let\startpage\relax
\let\finishpage\relax\let\publishdate\relax\let\receiveddate\relax
\let\reviseddate\relax\let\accepteddate\relax\let\theasciititle\relax
\let\theasciiauthors\relax\let\theasciiaddress\relax
\let\theasciiabstract\relax

\let\theerratum\relax\let\theasciiemail\relax
\let\theshortauthors\relax\let\theshorttitle\relax

\def\startpage{1}\def\finishpage{15}\def\thepapernumber{77}

\volumenumber{2}
\volumename{Proceedings of the Kirbyfest}
\volumeyear{1999}

\long\def\maketitlep{   

\count0=\startpage

\gtm\nl        
{\small Volume \thevolumenumber: \thevolumename\nl 
\ifx\theerratum\relax\else Erratum \erratumnumber\nl\fi
Pages \startpage--\finishpage\nl}

\vglue 0.1truein   

{\parskip=0pt\leftskip 0pt plus 1fil\def\\{\par\smallskip}{\ifplaintex\large
\else\Large\fi\bf\thetitle}\par\medskip}   
\vglue 0.05truein 

{\parskip=0pt\leftskip 0pt plus 1fil\def\\{\par}{\sc\theauthors}
\par\medskip}%
 
\vglue 0.03truein 

{\small\leftskip 25pt\rightskip 25pt{\bf Abstract}\stdspace\theabstract

{\bf AMS Classification}\stdspace\theprimaryclass
\ifx\thesecondaryclass\relax\else; \thesecondaryclass\fi\par
{\bf Keywords}\stdspace \thekeywords\par}\vglue 7pt

}   

\font\phead=cmsl9 scaled 950
\font=cmsl9 scaled 1050
\font\pnum=cmbx10 scaled 913
\font\lnum=cmbx10 
\font\pfoot=cmsl9 scaled 950
\font=cmsl9 scaled 1050
\ifplaintex
\headline{\vbox to 0pt{\vskip -4.5mm\line{\small\phead\ifnum
\count0=\startpage ISSN 1464-8997 (on line)
1464-8989 (printed) \hfill {\pnum\folio}\else\ifodd\count0\def\\{ }%
\ifx\theshorttitle\relax\thetitle\else\theshorttitle\fi\hfill{\pnum\folio}
\else\def\\{ and }{\pnum\folio}\hfill\ifx\theshortauthors\relax\theauthors
\else\theshortauthors\fi\fi\fi}\vss}}
\footline{\vbox to 0pt{\vglue 0mm\line{\small\pfoot\ifnum\count0=\startpage
Published \publishdate:\qua\copyright\ \gtp\hfill\else
\gtm, Volume \thevolumenumber\ (\thevolumeyear)\hfill\fi}\vss
}}
\else
\makeatletter
\def\@oddhead{{\small\lhead\ifnum\count0=\startpage ISSN 1464-8997 (on line)
1464-8989 (printed) \hfill {\lnum\number\count0}\else\ifodd\count0
\def\\{ }\ifx\theshorttitle\relax \thetitle \else\theshorttitle\fi\hfill
{\lnum\number\count0}\else\def\\{ and }{\lnum\number\count0}
\hfill\ifx\theshortauthors\relax 
\theauthors\else\theshortauthors\fi\fi\fi}}\def\@evenhead{@oddhead}
\def\@oddfoot{\small\lfoot\ifnum\count0=\startpage Published \publishdate:\qua\copyright\ \gtp\hfill\else
\gtm, Volume \thevolumenumber\ (\thevolumeyear)\hfill\fi}
\def\@evenfoot{@oddfoot}
\makeatother
\fi

\let\maketitlepage\maketitlep
\let\makeshorttitle\maketitlepage
\let\maketitle\maketitlepage

\newwrite\gtoutfile
\long\gdef\makeheadfile{  
{\def\\{, }\def\s{ }
\immediate\openout\gtoutfile head.xxx
\immediate\write\gtoutfile{To: math@arxiv.org}
\immediate\write\gtoutfile{Subject: put OR rep NNNNN:ppppp}
\immediate\write\gtoutfile{--text follows this line--}
\immediate\write\gtoutfile{Proxy-for: \ifx\theasciiauthors\relax
\theauthors\else\theasciiauthors\fi\s<\ifx\theasciiemail\relax\theemail\else\theasciiemail\fi>}
\immediate\write\gtoutfile{\noexpand\\}
\immediate\write\gtoutfile{Authors: \ifx\theasciiauthors\relax
\theauthors\else\theasciiauthors\fi}
{\def\\{ }\immediate\write\gtoutfile{Title: \ifx\theasciititle\relax
\thetitle\else\theasciititle\fi}}
\immediate\write\gtoutfile{Subj-class: GT or SG, GR etc}
\immediate\write\gtoutfile{MSC-class: \theprimaryclass\ifx\thesecondaryclass\relax\else, \thesecondaryclass\fi}
\immediate\write\gtoutfile{Journal-ref: Geom. Topol. Monogr. \thevolumenumber\s
(\thevolumeyear) \startpage-\finishpage}
\immediate\write\gtoutfile{Comments: Published by Geometry and Topology Monographs at}
\immediate\write\gtoutfile{\s\s\s  http://www.maths.warwick.ac.uk/gt/GTMon\thevolumenumber/paper\thepapernumber.abs.html}
\immediate\write\gtoutfile{\noexpand\\}
\immediate\write\gtoutfile{}
\ifx\theasciiabstract\relax
\immediate\write\gtoutfile{\theabstract}\else
\immediate\write\gtoutfile{\theasciiabstract}\fi
\immediate\write\gtoutfile{}
\immediate\write\gtoutfile{\noexpand\\}
\immediate\write\gtoutfile{}
\immediate\closeout\gtoutfile}}  

\def\maketitlepage{\maketitlep\makeheadfile}
\let\makeshorttitle\maketitlepage
\let\maketitle\maketitlepage

\volumenumber{4}
\volumename{Invariants of knots and 3-manifolds (Kyoto 2001)}
\volumeyear{2002}
\papernumber{14}
\pagenumbers{215}{233}
\received{7 April 2002}
\revised{12 October 2002}
\accepted{10 September 2002}
\published{13 October 2002}
 
\usepackage{amsmath,amssymb,amscd,epsf}
\def\R{\mathbb R}

\def\Ge{\varepsilon}

\def\Gs{\sigma}

\newcommand{\cY}{\mathcal Y}

\newtheorem{thm}{Theorem}[section]
\newtheorem{lem}[thm]{Lemma}

\theoremstyle{remark}
\newtheorem{rem}[thm]{Remark}
\theoremstyle{definition}
\newtheorem{defn}[thm]{Definition}
\newtheorem{ex}[thm]{Example}

\def\sminus{\smallsetminus}
\def\dd{\partial}

\def\cd{chord diagram}
\def\eq{equivalence}
\def\eqn{\underset{n}{\sim}}

\begin{document}

\title{Cubic complexes and finite type invariants}

\author{Sergei Matveev\\Michael Polyak}

\address{Department of Mathematics, Chelyabinsk State
University\\Chelyabinsk, 454021, Russia\\{\rm and}\\Department of Mathematics, 
Technion - Israel Institute of
Technology\\32000, Haifa, Israel}

\asciiaddress{Department of mathematics, Chelyabinsk State
University\\Chelyabinsk, 454021, Russia\\and\\Department of Mathematics, 
Technion - Israel Institute of
Technology\\32000, Haifa, Israel}

\email{matveev@csu.ru, polyakm@math.technion.ac.il}

\begin{abstract}
Cubic complexes appear in the theory of finite type invariants so
often that one can ascribe them to basic notions of the theory.  In
this paper we begin the exposition of finite type invariants from the
`cubic' point of view.  Finite type invariants of knots and homology
3-spheres fit perfectly into this conception.  In particular, we get a
natural explanation why they behave like polynomials.
\end{abstract}



\maketitle

\section{Introduction}

Polynomial functions play a fundamental role in mathematics.
While they are usually defined on Euclidean spaces, linear and even
quadratic maps are commonly considered for more general spaces, for
example for abelian groups.
This observation leads to a natural question: on which spaces can one
define polynomial functions, and which structure is required for that?
Certain hints pointing to a possible answer can be extracted from the
theory of difference schemes on cubic lattices.
For example, a continuous function is linear, if its forward second
difference derivative at any point $x_0$ vanishes, i.e.\
$f(x_0+x_1+x_2)-f(x_0+x_1)-f(x_0+x_2)+f(x_0)=0$ for any $x_1$ and $x_2$.
Similarly, quadratic functions are characterized by the identity
$f(x_0+x_1+x_2+x_3)-f(x_0+x_1+x_2)-f(x_0+x_1+x_3)-f(x_0+x_2+x_3)+
f(x_0+x_1)+f(x_0+x_2)+f(x_0+x_3)-f(x_0)\equiv0$.
It should be clear now how to generalize this to higher degrees:
\begin{thm}
\label{eucl}
A continuous function $f:\R^d\to\R$ is polynomial of degree less
than $n$ if and only if $\sum_\Gs (-1)^{|\Gs|} f(x_\Gs)=0$ for any
$x_0,x_1,\dots,x_{n}\in\R^d$. Here the summation is over all
$\Gs=(\Gs_1,\dots\Gs_{n+1})\in\{0,1\}^{n}$, $|\Gs|=\sum_i\Gs_i$
and $x_\Gs=x_0+\sum_i\Gs_i x_i$.
\end{thm}

Note that $\Gs$ are nothing more than the vertices of the standard
$n$-cube $[0,1]^n$ and $x_\Gs$ are their images under an affine map
$\varphi:[0,1]^n\to\R^d$.
Alternatively, if we replace the cube $[0,1]^n$ by $[-1,1]^n$, we
get the following characterization of polynomial functions:
$$\sum_\Gs\Gs_1\dots\Gs_n f(x_\Gs)\equiv0,$$
where $\Gs=(\Gs_1,\dots\Gs_n)\in\{-1,1\}^{n}$ and
$x_\Gs=x_0+\sum_i\Gs_i x_i$.
This formula corresponds to vanishing of $n$-th central difference
derivatives of $f$ at $x_0$.

Both formulas have the same meaning: the function is polynomial of
degree less than $n$, if the alternating sum of its values on the
vertices of every affine (possibly degenerate) $n$-cube in $\R^d$
vanishes.
Therefore, one may expect that given a set $X_n$ of "$n$-cubes" in
a space $W$, we may define a notion of a polynomial function $W\to\R$.
What are such $n$-cubes and how should they be related for different $n$?
An appropriate object is well-known in topology under the name of
a cubic complex.

While cubic complexes were used in topology for decades, their
relation to polynomials became apparent only recently in the
framework of finite type invariants. It turns out that cubic
complexes underlie so many properties of finite type invariants,
that one may ascribe them to basic notions of the theory. Probably
M.~Goussarov was one of the first to notice this
relation and realize its importance in full generality; the second
author learned this idea from him in 1996. This relation was also
noticed and discussed in an interesting unpublished preprint
\cite{FRS}. In this paper we begin the exposition of the theory of
finite type invariants from the ``cubic'' point of view.
Finite type invariants of knots and homology 3-spheres fit perfectly
into this conception.

\rk{Acknowledgements}The first author is partially supported by grants
E00-1.0-2.11, RFBR-02-01-01-013, and UR.04.01.033.  The second author
is partially supported by the Israeli Science Foundation, grant
86/01. The final version of the paper was written when both authors
visited the Max-Planck-Institut f\"ur Mathematik in Bonn.

\section{Cubic complexes}

\subsection{Semicubic complexes}
Simplicial complexes, i.e., unions of simplices in $\R^d$, are
widely used in topology.
While somewhat less intuitively clear, a notion of a semisimplicial
complex (see \cite{ssc}) is used in situations when the number of
simplices is infinite, especially locally.
Semicubic complexes are similar to semisimpli\-cial ones.
The only difference is that instead of simplices we take cubes.

\begin{defn} A {\em semicubic complex} $\bar X$ is a sequence of
arbitrary sets and maps $\ldots\Rightarrow X_n\Rightarrow
X_{n-1}\Rightarrow\ldots\Rightarrow X_0$. Here each arrow $X_n
\Rightarrow  X_{n-1}$ stands for $2n$ maps $\dd^{\Ge}_i \colon
X_n\to X_{n-1}$, $1\leq i \leq n$, $\Ge=\pm$, called the {\em
boundary operators}. The boundary operators are required to
commute after reordering: if $i>j$, then
$\dd^{\Ge_1}_j\dd^{\Ge_2}_i=\dd^{\Ge_2}_{i-1}\dd^{\Ge_1}_j.$
Elements of $X_n$ are called {\em $n$-dimensional cubes}, the maps
$\dd^{\pm}_i$ are {\em boundary operators}. A pair $\dd^-_i(x)$,
$\dd^+_i(x)\in X_{n-1}$ form the $i$-th pair of the {\em opposite
faces} of $x\in X_n$. One can consider also   a semicubic complex
$\bar X$ as a {\em semicubic structure} on the set $X_0$.
\end{defn}

\begin{figure}[ht!]
\centerline{\epsfysize3cm\epsfbox{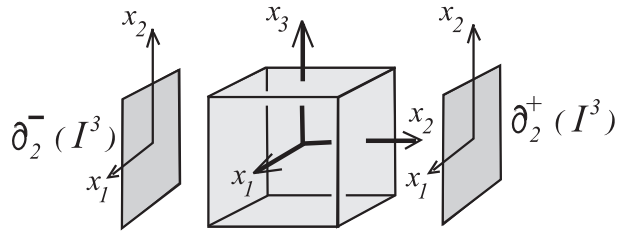}}
\caption{Boundary operators} \label{bouope}
\end{figure}

The above relations between maps mimic the usual identities for
the standard cube $I^n=\{(x_1,\ldots,x_n)\in\R^n\colon -1\leq
x_i\leq 1\}$ with vertices $(\pm 1,\pm 1,\ldots,\pm 1)$. Each time
when we take an $(n-1)$-face, we identify it with the standard
cube of dimension $n-1$ by renumbering the coordinate axes
monotonicaly. See Fig.~\ref{bouope}.

It follows from the commutation relations that any superposition
of $n$ boundary operators taking $X_n$ to $X_0$ coincides with a
monotone superposition $\dd^{\Ge_1}_1\dd^{\Ge_2}_2$ $\ldots
\dd^{\Ge_n}_n$ (or with a superposition
$\dd^{\Ge_n}_1\dd^{\Ge_{n-1}}_1\ldots \dd^{\Ge_1}_1$, whichever
you like). Therefore, any cube $x\in X_n$ has $2^n$ 0-dimensional
vertices, if we count them with multiplicities.   The set of all
vertices is naturally partitioned into two groups: we set a vertex
$\dd^{\Ge_1}_1\dd^{\Ge_2}_2\ldots\dd^{\Ge_n}_n$ to be {\em
positive} if $\Ge_1 \Ge_2 \ldots \Ge_n$ is $+$, and {\em negative}
otherwise.

By a map of a semicubic complex $X$ to a semicubic complex $Y$ we
mean a sequence $\bar \psi =\{ \psi_n \} $  of maps $\psi_n\colon
X_n\to Y_n, 0\leq n< \infty$, such that they commute with the
boundary operators, i.e., $\partial _i^\varepsilon \psi_n=
\psi_{n-1}\partial _i^\varepsilon$ for all $n\geq i\geq 1$ and
$\varepsilon=\pm $.

Evidently, semicubic complexes and maps between them form a
category.

\subsection{Incidence complexes}
To exclude the ordering of boundary operators (which is intrinsic
to semicubic complexes but often is inessential), we define cubic
complexes in more general terms  of {\em incidence relations}.

Let $A,B$ be an ordered pair of arbitrary sets. By an {\em
incidence relation} between $A$ and $B$ we mean any subset $R$ of
$A\times B$. If $(a,b)\in R$ for  some $a\in A, b\in B$, then
  we write $a\succ b$. The same   notation $A\succ B$ will be used
for indicating that $A,B$ are equipped with a fixed incidence
relation.
\begin{defn}
An {\em incidence complex }  ${\bar X }$ is a sequence $\dots\succ
X_n \succ X_{n-1}\succ \ldots\succ X_0 $ of arbitrary sets and
incidence relations between neighboring sets. Elements of $X_n$
are called {\em $n$-dimensional cells}. A cell $c_m\in X_m$ is a
{\em face} of a cell $c_n\in X_n, 0\leq m\leq n,$ if there exist
cells $c_{n-1}, \ldots , c_{m+1}$ such that $c_n\succ c_{n-1}\succ
\ldots \succ c_{m }$.
\end{defn}

By a map $\bar \psi \colon  {\bar X }  \to \bar  X'$ between two
incidence complexes  we mean a sequence of maps $\psi_n\colon
X_n\to X'_{n}$ such that $c_m\in X_m$ is a face of $c_n\in X_n$
implies $\psi_m(c_m)\in X'_m$ is a face of $\psi_n(c_n)\in X'_n$.

\begin{ex}
 The standard cube $I^n$ and all its faces form an incidence complex
${\bar I}^n$ in an evident manner. Obviously, each face $ I^m$ of
  $I^n $ is a cube     such that the inclusion  $I^m\subset I^n $
  induces an inclusion of the corresponding incidence complexes.
\end{ex}

\begin{defn} Let ${\bar X }$ be an incidence complex. Then any map
$\bar \varphi \colon {\bar  I}^n\to {\bar X }$  (considered as a
map between incidence complexes) is called a {\em cubic chart} for
$\bar X$.
\end{defn}
Evidently, for any face $ I^m$ of $  I^n, 0\leq m\leq n,$ the
restriction $\bar \varphi |_{\bar I^m}$  of any cubic chart $\bar
\varphi \colon {\bar I}^n\to {\bar X }$   is also a cubic chart.

\begin{defn} An incidence complex  ${\bar X }$ equipped with a set $\Phi$ of cubic charts
  for $\bar X$   is called   {\em cubic  } if the following holds:
\begin{enumerate}
\item For any cube $x\in X_n$ there is at least one chart
$\bar \varphi  \colon {\bar  I}^n\to \bar X$ in $\Phi  $ such that
$\varphi_n ( I^n) =x$. We will say that $\varphi $ {\em covers}
$x$.
\item The restriction of any cubic chart $\bar \varphi  \colon {\bar  I}^n\to
\bar X$ in $\Phi$ onto any face subcomplex $\bar  I^m$ of $\bar
 I^n$ belongs to $\Phi$.
\item For any two charts $\bar \varphi_1\colon {\bar  I }^n\to \bar X$ in $\Phi $,
$\bar \varphi_2\colon {\bar  I}^n\to \bar X$ which cover the same
cube $x\in X_n$ there exists   a combinatorial isomorphism (called
{\em a transient map}) $\bar \psi \colon  \bar
 I^n \to  \bar  I^n$ such that $\bar \varphi_1=\bar
\varphi_2 \bar \psi$.
\end{enumerate}
\end{defn}

\subsection{Oriented cubic complexes}
Of course, oriented cubic complexes are composed from oriented
cubes, but our definition of orientation of a cube drastically
differs from the usual one.

\begin{defn} An {\em orientation} of a cube   is an orientation
of all its edges such that all the parallel edges of the cube are
  oriented coherently. It means that if two edges are related by
a parallel translation, then   so are their orientations.
\end{defn}

To  any vertex $v$ of an oriented $n$-dimensional cube one can
assign  a sign $+$ or $-$, depending on the number of edges
outgoing of $v$: \ $+$ if it is even, and $-$ if odd.
Also, every pair of opposite $(n-1)$-dimensional faces of the
cube consists of a {\em negative} and a {\em positive} face.
We distinguish them by the behaviour of orthogonal edges: an
$(n-1)$-dimensional face is negative if all the orthogonal
edges go out of its vertices.
If all the orthogonal edges are incoming, then the face is positive.

\begin{figure}[ht!]
\centerline{\epsfysize 3.6cm\epsfbox{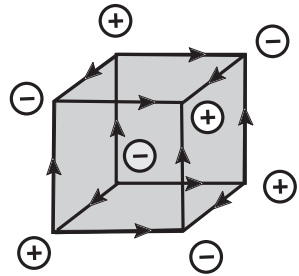}} \caption{The
signs of vertices. They will be used for taking alternative sums.} 
\label{signs}
\end{figure}

 Note that the  standard cube
$ I^n$ is equipped with the canonical orientation induced by the
orientations of
 the coordinate axes. See Fig.~\ref{signs} for the
signs of its vertices.

\begin{rem}\label{obv} Obviously, any orientation-preserving
isomorphism $I^n\to    I^n$
preserves the signs of vertices and
$(n-1)$-dimensional faces, and keeps fixed the {\em source vertex}
(having only outgoing edges) as well as the {\em sink vertex}
(having only incoming ones).
\end{rem}

\begin{defn} A cubic  incidence complex $\bar X$ is {\em oriented},
if every two its charts covering the same $n$-cube are related by
an orientation-preserving transient map $\bar I^n \to \bar I^n$.
\end{defn}

It follows from Remark~\ref{obv} that if a cubic incidence
complex $\bar X$ is oriented, then all vertices and
$(n-1)$-dimensional faces of any cube $x\in X_n$ have correctly
defined signs. Of course, if  $x$ has less than $2^n$ vertices,
then some of them  have several signs. Source and sink vertices of
$x$ are also defined, as well as its positive and negative faces.

 We say that a map $\bar \psi \colon \bar X \to \bar X'$  between
two oriented cubic incidence complexes is orientation-preserving,
if for any cube $x\in X_n$ of $\bar  X$ there exist cubic charts
$\bar \varphi \colon \bar I^n \to \bar X$ of $\bar X$ and  $\bar
\varphi' \colon \bar I^n \to \bar X'$ of $\bar X'$ and an
orientation-preserving map $\bar \psi' \colon \bar I^n \to \bar
I^n$ such that  $\bar \varphi$ covers $x$ and  $\bar \psi \bar
\varphi=\bar \varphi' \bar \psi'$. Of course, oriented cubic
complexes and orientation-preserving maps form a category.

\begin{rem} \label{semi-cub} Any semicubic complex $\bar Y$
determines an oriented cubic complex $\bar X$: we simply forget
about ordering of the boundary operators, preserving the information
on positive and negative $(n-1)$-dimensional faces. Vice versa,
any oriented cubic complex $\bar X$ determines a semicubic complex
$\bar Y$ as follows. The set $Y_n$ consists of all cubic maps
$\bar I^n \to \bar X$ which are related to cubic charts of $\bar X$
by orientation-preserving transient maps. The boundary operators
$\partial^{\pm }_i$ are defined by taking restrictions onto
positive and negative $i$-th faces of $I^n$. Both constructions
are functorial.
\end{rem}

\begin{ex} Let $W$ be a topological space.
Then {\em singular cubes} in $W$, i.e., continuous maps $f\colon
[-1,1]^n\to W$ of standard cubes into $W$, can be organized into a
semicubic complex as well as into  an oriented cubic complex $\bar
X=\bar X(W)$ in an evident way: $X_n$ is the set of all singular
cubes of dimension $n$, and the boundary operators, respectively,
incidence relations are given by taking restrictions onto the
faces.
\end{ex}

Other examples are discussed in Section \ref{sec_ex}. As we have
seen in Remark~\ref{semi-cub},  semicubic complexes and oriented
cubic incidence complexes are related very closely.  Further on
 we will use the semicubic complexes, but occasionally return to
 oriented  ones. For brevity, in both cases we will call them
``cubic complexes''.

\subsection{Cubes vs simplices}
Cubes     enjoy all good properties of  simplices and have the
following advantages:
\begin{enumerate}
\item Each face of a cube has the opposite face;
\item Two cubes with a common face can be glued together into a
new cube (well, parallelepiped, but it does not matter);
\item The direct product of two cubes of dimensions $m$ and $n$
is a cube of dimension $m+n$.
\end{enumerate}

The above properties of cubes may be included as axioms.
We will say that a  cubic complex   is {\em good}, if
\begin{enumerate}

\item For each $n$ and $i, 1\le i\le n,$
 there is an
involution $J_i:X_n\to X_n$, such that
 $\dd^\Ge_j J_i=J_{i-1}\dd^\Ge_j$ for $j<i$,
$\dd^\Ge_j J_i=J_{i}\dd^\Ge_j$ for $j>i$, and
 $\dd^\Ge_i J_i=\dd^{-\Ge}_i$;

\item For each $x,y\in X_n$ with $\dd^+_i(x)=\dd^-_i(y)$
there is $x\circ y\in X_n$ such that $\dd^+_i(x\circ
y)=\dd^+_i(y)$ and $\dd^-_i(x\circ y)=\dd^-_i(x)$;

\item For each $x\in X_n$, $y\in X_m$ with $\dd^-_1\ldots\dd^-_n(x)
=\dd^-_1\ldots\dd^-_m(y)$ there is $xy\in X_{n+m}$ such that
$\dd^-_1\ldots\dd^-_n(xy)=y$ and $\dd^-_{n+1}\ldots\dd^-_{n+m}(xy)=x$.
\end{enumerate}

These axioms guarantee a rich algebraic structure (duality,
composition, and product) on good cubic complexes.  One can easily
show that any cubic complex can be embedded into a good cubic
complex. However, in all interesting examples which we presently
know the cubic complexes are good. Note that axioms 2 and 3
descend to the level of oriented cubic complexes.

Another important advantage of cubes over simplices, as we
will see below, is that they turn out to be extremely useful
for a study of polynomial functions.

\section{Finite type functions and $n$-equivalence}

\subsection{Functions on vertices}
Just as for semisimplicial complexes, for a semicubic complex
$\bar X$ one may consider {\em $n$-chains} $C_n(\bar X)$, i.e.,
linear combinations of $n$-cubes with, say, rational coefficients.
Any function $f$ on $X_n$ extends to a function on $n$-chains
$C_n(\bar X)$ by linearity. The boundary operators
$\dd^\Ge_i:X_n\to X_{n-1}$, picked with an appropriate signs, may
be combined into a differential $\sum_i(-1)^i(\dd^+_i-\dd^-_i)$ on
chains. This differential brings us to homology groups.

Having in mind a study of polynomial functions we, however, choose
the signs differently. Namely, we define the operator $\partial
\colon C_n(\bar X)\to C_{n-1}(\bar X)$ by
$$\dd=\frac{1}{n}\sum_i(\dd^+_i-\dd^-_i).$$ This operator does not
satisfy $\dd^2=0$; of a special interest for us will be $0$-chains
$\dd^n(x)$, for any $x\in X_n$. The chain $\dd^n(x)$ contains all
$2^n$ vertices of the $n$-cube $x$ with signs shown in
Fig.~\ref{signs}: $$
\dd^n(x)=\sum_{\Ge_1,\ldots,\Ge_{n}}\Ge_1\ldots\Ge_{n}
 \dd^{\Ge_1}_1\dd^{\Ge_2}_2\ldots\dd^{\Ge_{n}}_{n}(x)
 .$$
\begin{rem} Note that the sum $\frac{1}{n}\sum_i(\dd^+_i-\dd^-_i)$ does
not depend on the order of the summands. It means that $\partial$
is determined for any oriented cubic  complex.
\end{rem}

\subsection{Polynomials in $\R^d$}\label{sub_R}
Let us start with the following simple example.

Let $X_n$ consist of all affine maps $\varphi \colon I^n\to\R^d$
and the boundary operators  $\dd^{\pm}_i$ assign to each $\varphi$
its restrictions to the faces $\{ x_i=\pm 1\}$. Then $X_0$ can be
identified with $R^d$. We would like to interpret polynomiality of
functions $f:\R^d\to\R$ in terms of their values on $\dd^n(x)$,
$x\in X_n$.
In these terms Theorem \ref{eucl} may be restated as follows:

\begin{thm}
A continuous function $f:\R^d\to\R$ is polynomial of degree less
than $n$, if and only if $f(\dd^n(x))=0$ for all $x\in X_{n}$.
\end{thm}

Recall that the cubes may be highly degenerate: to obtain a
polynomial dependence on $i$-th coordinate, we consider affine
maps with the image of $I^n$ contained in a line parallel to
the $i$-th coordinate axis.

\subsection{Finite type functions}
The example above motivates the following definition, which works
for both semicubic and oriented cubic complexes.

\begin{defn}
Let $\bar X$ be a cubic complex, and $A$ an abelian group. A
function $f:X_0\to A$  is of {\em finite type of degree less than
$n$}, if for all $x\in X_n$ we have $f(\dd^n(x))=0$.
\end{defn}

\begin{rem}
This definition looks more familiar in terms of the dual {\em
cochain complex}
 $C^*(\bar X)$ of linear functions on $C_*(\bar X)$
with the coboundary operator $df(x)=f(\dd x)$ dual to $\dd$:
 a function $f\in C^0({\bar X})$ is of degree less than $n$ if $d^n(f)=0$.
\end{rem}

Note that   the cubic complex
$\bar U$ whose $n$-cubes are   affine functions $f\colon I^n\to R$
has the following   interesting property: any   finite type function on
$U_0$ is a polynomial.
This observation explains why in
many respects finite type functions  behave analogously to
polynomials \cite{BN2}.

\begin{figure}[ht!]
\centerline{\epsfysize 3cm\epsfbox{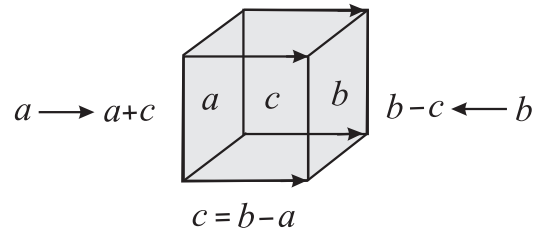}}
\caption{Jumping in positive direction we gain $c$, jumping back
we loose $c$. } \label{jump}
\end{figure}

Given a function $f:X_n\to A$ on $n$-cubes, sometimes one may
extend it to all the cubes of dimension $<n$, including vertices,
by the following descending method (first described in different
terms by Vassiliev \cite{Va} for a cubic complex  of knots).   By
a {\em jump} through an $n$-cube $x\in X_n$ we mean the transition
from an $(n-1)$-face of $x$ to the opposite one. The value
$c=f(x)$ is the {\em price} of the jump. We add $c$, if we jump
from a negative face to the opposite one, and subtract $c$, if the
jump is from a positive face to the opposite negative one.
See Fig.~\ref{jump}.

\begin{thm}
A function $f\colon X_n\to A$ extends to $X_{n-1}$ if and only if
the algebraic sum of the prices of any cyclic sequence of jumps is
zero.
\end{thm}

\begin{figure}[ht!]
\centerline{\epsfysize3cm\epsfbox{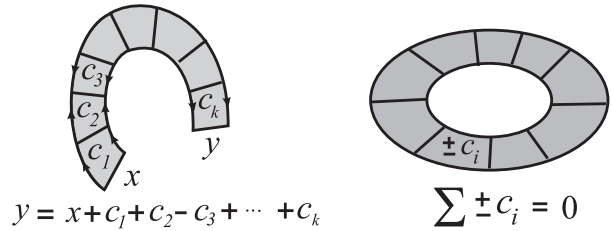}}
\caption{Prices of jumps and the cyclic condition }
 \label{desce}
\end{figure}

\begin{proof}
Call two $(n-1)$-cubes {\em parallel}, if one can pass from one to
the other by jumps. In each equivalence class we choose a
representative $r$, assign a variable, say, $y$ to it, and set
$f(r)=y$. Then we calculate the value of $f$ for any other cube
from the same equivalence class by paying $\pm f(x)$ for jumping
across any cube $x\in X_n$. It is clear that we get a correctly
defined function on $X_{n-1}$ if and only if the above cyclic
condition holds, see Fig.~\ref{desce}.
\end{proof}

One can look at the descending process as follows. To construct a
function of degree n, we start with the zero function $X_{n+1}\to
A$ and try to descend it successively to functions on $X_n ,
 X_{n-1} ,\dots, X_0$. At each step we create a lot of new
variables, and at each next step subject them to some linear
homogeneous restrictions. If the system has a nonempty solution
space, then we can descend further.

\begin{rem} The number of equations at each step
can be infinite.  Sometimes it can  be made finite by the
following two tricks. First, it suffices to consider only basic
cycles, which generate all cyclic sequences of $n$-cubes. Second,
if one cyclic sequences of $n$-cubes consists of  some
 faces of    a cyclic sequence of $(n+1)$-cubes, then the
sequence of the opposite   faces is also cyclic and  has the same
sum of jumps. Therefore, it suffices to consider only one of these
two  chains.

 The question when the solution space is
always nonempty (i.e., when the descending process gives us
nontrivial invariants) is usually very hard. For the case of
real-valued finite type invariants of knots (see the next
section), when the number  of variables and equations at each step
can be made finite,  the affirmative answer follows from the
Kontsevich theorem \cite{Kon}. For the case of finite type
invariants of homology spheres the number of variables is
infinite, which makes this case especially difficult. See
\cite{GGP}, where the authors managed to get rid off all but
finitely many variables by borrowing additional relations from the
lower levels.
\end{rem}

\subsection{N-\eq\ and \cd s}
For any cubic complex  one may define a useful notion   of $n$-\eq :

\begin{defn}
Let $\bar X$ be a cubic complex. Elements $x,y\in X_k$ are {\em
$n$-equivalent}, if there exists an $(n+k+1)$-chain $z\in C_{n+k+1}(\bar X)$,
such that $\dd^{n+1}(z)=y-x$. We denote $x\eqn y$.
\end{defn}

\begin{rem}
Our definition is somewhat different from the one used by\break
M.~Goussarov for links and 3-manifolds. He defines the notion of
$n$-\eq\ only on $X_0$ and uses a certain additional geometrical
structure present in these cubic complexes. Roughly speaking, his
relation is generated by $(n+1)$-cubes $z$ with
$\dd^{n+1}(z)=(-1)^{n+1}(y-x)$, but only of a special type, namely
such that $\dd^+_1\dots \dd^+_{n+1}(z)=y$
and $\dd^{\Ge_1}_1\dots \dd^{\Ge_{n+1}}_{n+1}(z)=x$ otherwise.
In other words, all the vertices of $z$ should coincide with an
0-cube $x$, except the unique sink vertex $y$.
One may show that $x$ and $y$ are $n$-equivalent in a
sense of Goussarov, if and only if there exist two $(n+1)$-cubes
$z_x$ and $z_y$ all vertices of which coincide, except for the sink
vertex, which is $x$ for $z_x$ and $y$ for $z_y$.

While {\em a priori} Goussarov's definition is finer, these
definitions are equivalent for knots; also, if the theorems
announced in \cite{Gu2, Gu3} are taken in the account, these
definitions should coincide for string links and homology
cylinders.
\end{rem}

It is easy to see that:

\begin{lem}
$n$-\eq\ is an \eq\ relation.
\end{lem}

\begin{rem}
Given a product on the set $X_0$, often the classes of
$n$-equival\-ence form a group. See \cite{Gu1, Gu2, Gu3}
for groups of knots, string links, and homology spheres.
\end{rem}

If a cubic complex  has more than one class of $0$-\eq , one may
also consider a restriction of the theory of finite type functions
to some fixed class of $0$-\eq . There exists also a more general
theory of partially defined finite type invariants, see
\cite{Gu2}.

The following simple theorem shows that functions of degree $\le n$
are constant on classes of $n$-\eq :

\begin{thm}
Let $\bar X$ be a cubic complex , $A$ an abelian group, and let
$f:X_0\to A$ be a function of degree $\le n$. Then for any $x,y\in
X_0$ such that $x\eqn y$ we have $f(x)=f(y)$.
\end{thm}

\begin{proof}
By the definition of $n$-\eq , there exists an $(n+1)$-chain
$z$ such that $\dd^{n+1}(z)=y-x$.
But $f$ is of degree $\le n$, hence $f(\dd^{n+1}(z))=0$.
It remains to notice that $f(x)-f(y)=f(x-y)=f(\dd^{n+1}(z))$.
\end{proof}

\begin{rem}
The opposite is not true: in general functions of finite type do
not distinguish classes of $n$-\eq , i.e.,  the equality
$f(x)=f(y)$ for any $f$ of degree $\le n$  does not imply that
$x\eqn y$ (see \cite{Gu2} for the case of the link cubic complex).

In some important cases, however, e.g., for the case of the knot
cubic complex , functions of finite type do distinguish classes of
$n$-\eq , see \cite{St, Gu2, Gu3}. Goussarov \cite{Gu2, Gu3} also
announced similar results for the cubic complexes of string links
and homology cylinders, but we do not know what were his ideas on
the subject and no proofs seem to be known.

A study of conditions under which functions of finite type distinguish
classes of $n$-\eq\ present an important problem.
\end{rem}

It is also interesting to consider the quotients of $n$-\eq\
classes by the relation of $(n+1)$-\eq .

A closely related space $H_n(\bar X)$ of {\em \cd s} of a cubic
complex\ $\bar X$ is defined as $H_n(\bar X)=C_n/\dd(C_{n+1})$,
where $C_n$ are $n$-chains in $\bar X$. The {\em weight system}
$\{f(h)|h\in H_n\}$ of a function $f$ of degree $n$ is defined by
setting $f(h)=f(\dd^n z)$ for any representative $z\in h$. For any
other representative $z+\dd z'\in h$, $z'\in X_{n+1}$ we have
$f(\dd^n(z+\dd z'))=
f(\dd^n(z))+f(\dd^{n+1}(\dd^{n+1}(z'))=f(\dd^n(z))$, since $f$ is
of degree $n$.

\section{Examples of cubic complexes}
\label{sec_ex}

The examples in this section mimic the following definition of an
$n$-dimensional cube $x$ in $\R^d$ with sides parallel to the axes.
Let $x$ be a sequence of $d$ symbols, $d-n$ of which are real numbers
and $n$ are $*$, together with $n$ pairs $(s_i^-,s_i^+)$ of real
numbers.
An $(n-1)$-dimensional face $\dd_i^\pm x$ is obtained by
plugging $s_i^\pm$ in $x$ instead of $i$-th $*$.
To obtain a cube with sides which are not parallel to the
axes, instead of $i$-th $*$ coordinate (varying from $s_i^-$ to
$s_i^+$ along $i$-th side), one may consider several $*$'s
united in $i$-th group.

\subsection{Cubic structure on a group}
Let $G$ be a group.
We define a cubic structure on $G$ in the following way.
An $n$-cube $x\in X_n$ is a word $g_0(a_1b_1)g_2\dots(a_nb_n)g_n$
which contains $n$ bracketed pairs $(a_i,b_i)\in G\times G$
separated by elements $g_0,\dots,g_n$ of $G$.
Elements of $X_0$ are identified with $G$.
The boundary operator $\dd^-_i$ changes the $i$-th bracket
$(a_i,b_i)$ into the product $a_ib_i$, while the boundary
operator $\dd^+_i$ changes the $i$-th bracket into $b_ia_i$.
Since each $\dd^+_i$ differs from $\dd^-_i$ by the transposition
of a pair of elements $a_i$ and $b_i$, it is easy to see that
$0$-\eq\ classes coincide with the abelinization $G/[G,G]$
of $G$.
More generally, $n$-\eq\ classes coincide with the double
cosets $G_n\backslash G/G_n$ of $G$ by the $(n+1)$-th lower
central subgroup $G_n$.
Here the lower central subgroups $G_n$ of a group $G$ are
generated by $[G,G_{n-1}]$, with $G_0=G$.

Another, closely related, but somewhat more general cubic structure
on $G$ may be defined as follows.
Set $X_n$ to be a free product of $n+1$ copies $G^0$, $G^1$,\dots,
$G^n$ of $G$ (with $X_0=G^0$ identified with $G$).
The boundary operator $\dd^-_i$ is the projection of $G^i$ to $1$.
The boundary operator $\dd^+_i$ is the map identifying $i$-th
copy $G^i$ of $G$ with $G^0$ (and renumbering the other copies
of $G$ by $j\to j-1$ for $j>i$).
It is easy to see that linear functions $f:G\to A$ such that $f(1)=0$
are exactly homomorphisms of $G$ into $A$ (since vanishing of $f$ on a
2-cube $gh$ with $g\in G^1$, $h\in G^2$ implies $f(ab)-f(a)-f(b)+f(1)=0$).

\begin{thm}
Classes of $n$-\eq\ of the above cubic structure on $G$
coincide with the cosets $G_n\backslash G/G_n$ of $G$ by the $n$-th
lower central subgroup $G_n$.
\end{thm}

\begin{proof}
Indeed, from any element of the form $y=x\cdot ghg^{-1}h^{-1}$
we may construct a 2-cube as follows.
Write $xghg^{-1}h^{-1}$ as an element $z$ of the free product
of three copies $G^0$, $G^1$ and $G^2$ of $G$, with $x\in G^0$,
$g,g^{-1}\in G^1$, and $h,h^{-1}\in G^2$.
An application of $\dd^-_1$ removes $g$ and $g^{-1}$, leaving
$xhh^{-1}=x$; similarly, an application of $\dd^-_2$ removes
$h$ and $h^{-1}$, leaving $xgg^{-1}=x$.
Thus $\dd^+_1\dd^+_2(xghg^{-1}h^{-1})=y$ and
$\dd^{\Ge_1}_1\dd^{\Ge_2}_2(xghg^{-1}h^{-1})=x$
otherwise, so the 2-cube $z$ satisfies $\dd^2(z)=y-x$.
Hence $x$ and $y$ are $1$-equivalent.

In a similar way, if $y=x\cdot c$ differs from $x$ by an $n$-th
commutator $c$, then we may write it as an element $z$ in a free
product of $n+2$ copies of $G$, such that again $\dd^-_i(z)=x$
for any $1\le i\le n+1$, and hence $\dd^{n+1}(z)=y-x$ and $x\eqn y$.
The opposite direction is rather similar.
\end{proof}

The algebra $H_n$ of chord diagrams is in this case a free product
of $n$ copies of $G/[G,G]$.

\subsection{Cubic structures on trees, operads, and graphs}
Let $P$ be an operad. A cubic structure on $P$ may be
introduced by plugging some fixed $s^\pm_i$ in some
$x\in P$.
We will illustrate this idea on an example of the rooted
tree operad.

\begin{figure}[ht!]
\centerline{\epsfxsize 7cm\epsfbox{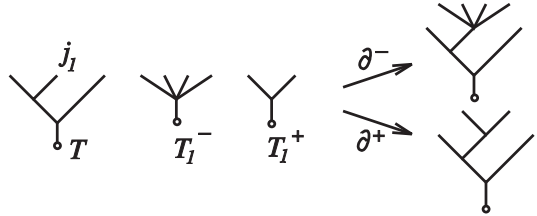}} \caption{Grafting
of trees} \label{trees}
\end{figure}

Define $x\in X_n$ as a collection $(T,T^\pm_i,j)$,
where $T$, $T_i^\pm$, $i=1,\dots,n$ are trees and
$j=j_1,\dots,j_n$ is a set of $n$ leaves of $T$.
The face $\dd^\pm_i x$ is obtained by attaching (grafting)
the root of $T^\pm_i$ to the $j_i$-th leaf of $T$.
See Figure \ref{trees}.

It would be interesting to study the theory of finite type
functions on this cubic complex . The example below shows that the
most basic functions of trees fit nicely in the theory of finite
type functions.

\begin{thm}
The number of edges and the number of vertices of a given valence
(as well as any of their functions) are degree one functions on
the cubic complex\ of trees. The number of $n$-leaved trees of
some fixed combinatorial type is a function of degree $\le n/2$.
\end{thm}

A similar cubic structure may be defined on graphs (using
insertions of some subgraphs $G^\pm_i$ in $n$ vertices). In
particular, using subgraphs which contain just one edge, we obtain
the following cubic structure. Define an $n$-cube $x\in X_n$ to be
a graph $G$ with $n$ marked vertices $v_1,\dots,v_n$, together
with two fixed partitions $s^+_i$, $s^-_i$ of edges incident to
$v_i$. The boundary operator $\dd^-_i x$ (resp. $\dd^+_i x$ acts
by inserting in $v_i$ a new edge, splitting it into two vertices
in accordance with the partition $s^-_i$ (resp. $s^+_i$) of the
edges. Here are some simple examples of finite type functions on
the cubic complex\ of graphs.

\begin{thm}
The number of edges, the number of vertices, and the number
of loops (as well as any of their functions) are degree
zero functions.
The number of vertices of some fixed valence is a degree
one function.
The number of edges with the endpoints being vertices of
some fixed valences is a degree two function.
The number of $n$-vertices subgraphs of some fixed
combinatorial type is a function of degree $\le n/2$.
\end{thm}

One of the relations in the algebra of \cd s for this cubic
complex\ is Stasheff's pentagon relation. We do  not know whether
there are any other relations. It may be also interesting to
investigate the relation of this cubic complex\ to the graph
cohomology.

\subsection{Vassiliev knot complex}\label{ex_knot1}
\begin{figure}[ht!]
\centerline{\epsfysize2cm\epsfbox{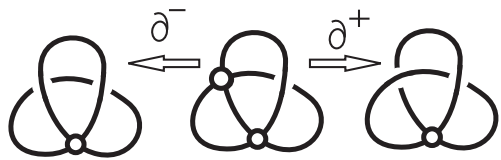}} \caption{A
2-singular knot and resolutions of a double point} \label{resol}
\end{figure}

Let $X_n$ consist of singular knots in $R^3$ having $n$
ordered transversal double points.
The  boundary operators  $\dd^{\pm}_i$ act by a positive,
respectively, negative resolution of $i$-th double point,
by shifting one string of the knot from the other.
See Fig.~\ref{resol}.
Here the resolution is positive, if the orientation of the fixed
string, the orientation of the moving string, and the direction
of the shift determine the positive orientation of $R^3$.
The vertices of an $n$-cube thus may be identified with $2^n$
knots, obtained from an $n$-singular knot by all the resolutions
of its double points.

Finite type functions for this complex are known as finite type
invariants of knots (also known as Vassiliev or
Vassiliev-Goussarov invariants), see \cite{BN1} for an elementary
introduction to the theory of Vassiliev invariants.

\subsection{More general knot complex}\label{ex_knot2}
\begin{figure}[ht!]
\centerline{\epsfxsize10cm\epsfbox{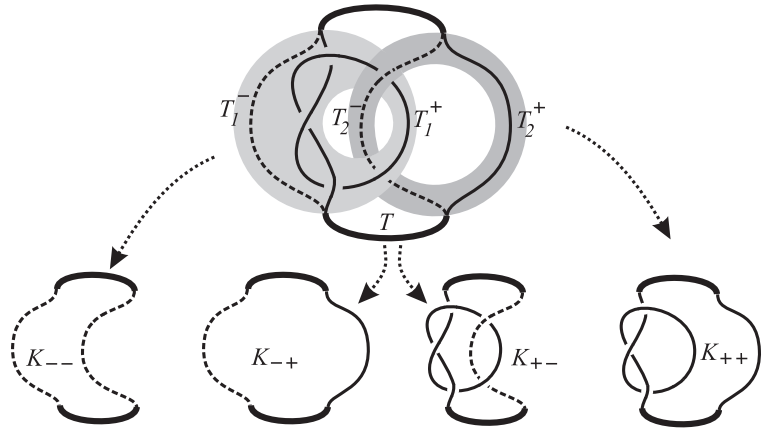}} \caption{A tangle
and its modifications} \label{modif}
\end{figure}

Let an element of $X_n$ be a knot together with a set of its $n$
fixed modifications.
More precisely, fix a set $H_1$, \dots, $H_n$ of disjoint
handlebodies in $\R^3$.
An $n$-cube is a tangle $T$ in $\R^3\sminus\cup H_i$, together with
a set of $2n$ tangles $T^\pm_i\subset H_i$, such that for any choice
of signs $\Ge_1,\ldots,\Ge_n$ the glued tangle
$K_{\Ge_1,\ldots,\Ge_n}=T\cup T^{\Ge_1}_1\cup\ldots\cup T^{\Ge_n}_n$
is a knot.
Here by a tangle in a manifold $M$ with boundary we mean a
1-dimensional manifold, properly embedded in $M$.
The boundary operators $\dd^\Ge_i$ act by forgetting $H_i$ and gluing
$T^\Ge_i$ to $T$.
See Figure \ref{modif}.
The vertices may be thus identified with $2^n$ knots
$K_{\Ge_1,\ldots,\Ge_n}$.

This cubic structure on knots was introduced by Goussarov
\cite{Gu4} under the name of "interdependent knot modifications".
From the construction (restricting the modifications to crossing
changes) it is clear that any finite type function in this theory
is a Vassiliev knot invariant. As shown in \cite{Gu4}, the
opposite is also true, so the finite type functions for this cubic
complex\ are exactly Vassiliev knot invariants (with a shifted
grading); see \cite{Gu4, BN3}.

It would be interesting to construct similar cubic complexes for
virtual knots and plane curves with cusps.

\subsection{Borromean surgery in 3-manifolds}
\begin{figure}
\centerline{\epsfysize4cm\epsfbox{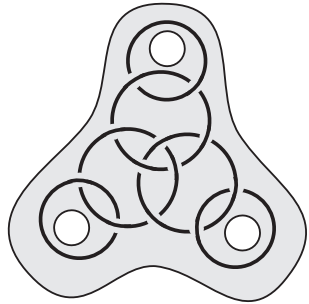}} \caption{A
$Y$-clasper} \label{borr}
\end{figure}

Let $H$ be a standard genus 3 handlebody presented as a 3-ball
with three index one handles attached to it.
Consider a 6-component link $L\subset H$ consisting of the Borromean
link $B$ in the ball and three circles which run along the handles
and are linked with the corresponding components of $B$,
see Fig.~\ref{borr}.
We equip $L$ with the zero framing.

\begin{defn}
An $n$-component  {\em $Y$-clasper} (or a $Y$-graph)
in a 3-manifold $M$ is a collection of $n$ embeddings $h_i:H\to M$,
$1\leq i\leq n$, such that the images $h_i(H)$ are disjoint.
\end{defn}

Let us construct a cubic complex\  as follows. $X_n$ consists of
all pairs $(M,\cY)$, where $M$ is a 3-manifold and $\cY$ an
$n$-component $Y$-clasper in $M$. The pairs are considered up to
homeomorphisms of pairs. The boundary operator   $\dd^{- }_i$ acts
by forgetting the $i$-th component $h_i$ of $\cY$ (the manifold
$M$ remains the same). The operator $\dd^{+ }_i$ also removes
$h_i$, but, in contrast to $\dd^{- }_i$, $M$ is replaced by the
new manifold obtained from $M$ by the surgery along $h_i(L)$. Such
a surgery is called {\em Borromean} (see \cite{Ma}). It is known
\cite{Ma} that one 3-manifold may be obtained from another by
Borromean surgeries (so belong to the same $0$-equivalence class)
if and only if they have the same homology and the linking pairing
in the homology. In particular, $M$ is a homology 3-sphere if and
only if it can be obtained from $S^3$ by Borromean surgeries. It
is easy to see that the sets $X_n$ together with operators
$\dd^{\pm}_i$ form a cubic complex .

Its finite type invariants are invariants of 3-manifolds in the
sense of \cite{Gu3, Ha}.
One may also restrict it to homology 3-spheres.

\subsection{Whitehead surgery in 3-manifolds}
There are several other approaches to the finite type invariants
of homology spheres.
They are based on surgery on algebraically split links \cite{Oh},
boundary links \cite{Ga}, blinks \cite{GL}, and so on.
All of them fit into the conception of cubic complexes and turn
out to be equivalent, see \cite{GGP}.

Here is a new approach, based on Whitehead surgery.

\begin{defn} An $n$-component $Y$-clasper $h_i:H\to M$ in a
3-manifold $M$ is a {\em Whitehead clasper}, if for each $i$ one
of the handles of $h_i(H)$ bounds a disc in $M\sminus \cup_j h_j(H)$
and the framing of this handle is $\pm1$.
\end{defn}

A surgery along a Whitehead clasper is called {\em Whitehead
surgery}; it was introduced in \cite{Ma} in different terms.
From the results of \cite{Ma} it follows that:

\begin{thm}
$M$ is a homology 3-sphere if and only if it can be obtained
from $S^3$ by surgery on a Whitehead clasper.
\end{thm}

We obtain the Whitehead cubic complex\  of homology 3-spheres by
considering only Whitehead $Y$-claspers in homology spheres in the
definition of the Borromean cubic complex\  above. It is easy to
see that the sets $X_n$ together with operators $\dd^{\pm}_i$ form
a cubic complex.

From the construction it is clear that all finite type invariants
of homology spheres of degree $<n$ in the sense of Borromean theory
above are also finite type invariants of degree $<n$ in the sense
of Whitehead surgery.
We expect the opposite to be also true (probably up to a degree
shift).
Considering this theory for arbitrary 3-manifolds, we get, however,
a theory which is finer than the theory based on the Borromean surgery.
The reason is that the Whitehead surgery preserves the triple cup
product in the homology, while the Borromean surgery in general
does not.
The study of this theory and its comparison with the theory
introduced in \cite{CM} seem to be promising.

\subsection{More on polynomiality}
In many examples (see above) there exist several different cubic
structures on the same space $X_0$.
However, in all presently known non-trivial examples the set of
finite type functions remains the same, up to a shift of grading.
See \cite{Gu4} for the case of knots, and \cite{GGP} for homology
3-spheres.

It would be quite interesting to understand better this "robustness"
of finite type functions, and to formulate conditions which would
imply such a uniqueness.

In conclusion we note that finite type invariants of knots and
homology spheres are obtained by the same schema as polynomials
(see Section \ref{sub_R}). This observation explains once again
their polynomial nature \cite{BN2}.
It is also worth noting a curious "secondary" polynomiality of finite
type invariants: any finite type invariant is a polynomial in
primitive finite type invariants.

Finally, let us remark that an oriented cubic complex (with n-cubes being
certain commutative diagrams of vector spaces) appear in the
construction of Khovanov's homology \cite{Kh} for the Jones polynomial. It
would be interesting to investigate Khovanov's construction from this
point of view.

\Addresses\recd

\end{document}